\documentclass[a4paper,12pt,spanish]{article}
\usepackage[T1]{fontenc}
\usepackage[latin1]{inputenc}
\usepackage{babel}
\usepackage{amsthm,amsmath}
\usepackage{amssymb}
\usepackage[all]{xypic}
\usepackage{enumerate}
\usepackage{a4wide}
\usepackage{hyperref}
\usepackage{latexsym}
\usepackage{enumerate}

\newtheorem{theorem}{Teorema}[section]
\newtheorem{proposition}[theorem]{Proposici\'{o}n}
\newtheorem{corollary}[theorem]{Corolario}
\newtheorem{lemma}[theorem]{Lema}
\newtheorem*{theorem*}{Teorema}
\newtheorem*{proposition*}{Proposici\'{o}n}
\newtheorem*{corollary*}{Corolario}
\newtheorem*{lemma*}{Lema}
\theoremstyle{definition}

\newtheorem{remark}[theorem]{Observaci\'{o}n}
\newtheorem*{remark*}{Observaci\'{o}n}

\newtheorem*{definition*}{Definici\'{o}n}

\newcommand{\cat}[1]{\mathcal{#1}}
\newcommand{\coring}[1]{\mathfrak{#1}}
\newcommand{\tensor}[1]{\otimes_{#1}}
\newcommand{\rcomod}[1]{\mathsf{Comod}_{#1}}
\newcommand{\rmod}[1]{\mathsf{Mod}_{#1}}

\newcommand{\rmodu}[1]{\overline{\mathsf{Mod}}_{#1}}

\newcommand{\cotensor}[1]{\square_{#1}}

\renewcommand{\hom}[3]{\mathrm{Hom}_{#1}(#2,#3)}

\newcommand{\rend}[2]{\mathrm{End}({#2}_{#1})}

\newcommand{\dostensor}[3]{#1 \tensor{#2} #3}
\newcommand{\trestensor}[5]{#1 \tensor{#2} #3 \tensor{#4} #5}
\newcommand{\fourtensor}[7]{#1 \tensor{#2} #3 \tensor{#4} #5 \tensor{#6} #7}
\newcommand{\abrir}[1]{e_{#1}\tensor{A}e^*_{#1}}
\newcommand{\coalg}[2]{{#1}_{#2}}

\begin{document}
\title{Com\'{o}nadas y coanillos de Galois\footnote{Investigaci{\'o}n realizada en el proyecto MTM2004-01406 <<M{\'e}todos algebraicos en Geometr{\'\i}a no conmutativa>> financiado por DGICYT y FEDER}}
\author{J. G\'omez-Torrecillas \\
\normalsize Departamento de \'{A}lgebra \\ \normalsize Universidad
de Granada\\ \normalsize E18071 Granada, Espa\~{n}a \\
\normalsize e-mail: \textsf{gomezj@ugr.es} }

\date{{\small Versi{\'o}n 1.2}}
\maketitle

\section*{Introducci\'{o}n}

La noci{\'o}n de coanillo fue introducida por M. E. Sweedler en
\cite{Sweedler:1975} con el objeto de formular y demostrar un
predual al Teorema de Jacobson-Bourbaki para extensiones de
anillos de divisi{\'o}n. Un argumento fundamental en
\cite{Sweedler:1975} es el siguiente: dados $E \subseteq A$
anillos de divisi{\'o}n, cada coideal $J$ del $A$--coanillo $A
\tensor{E} A$ da lugar a un coanillo cociente $\coring{C} =
A\tensor{E} A / J$. Si $g \in \coring{C}$ denota el elemento
<<como de grupo>> $1 \tensor{E} 1 + J$, entonces $D = \{ a \in A :
ag = ga \}$ es un anillo de divisi{\'o}n intermedio $E \subseteq D
\subseteq A$. Adem{\'a}s, se tiene el homomorfismo can{\'o}nico de
$A$--coanillos $\zeta : A \tensor{D} A \rightarrow \coring{C}$ que
lleva $1 \tensor{D} 1$ en $g$. Resulta de \cite[2.2 Fundamental
Lemma]{Sweedler:1975} que $\zeta$ es un isomorfismo de
$A$--coanillos, lo cual es, a la postre, b{\'a}sico para establecer la
correspondencia entre coideales de $A \tensor{E} A$ y extensiones
intermedias $E \subseteq D \subseteq A$ \cite[Fundamental
Theorem]{Sweedler:1975}. El citado Lema Fundamental de Sweedler
puede reemplazarse por el hecho de que el coanillo $A \tensor{D}
A$ resulta ser simple cosemisimple \cite[Theorem
4.4]{ElKaoutit/Gomez/Lobillo:2004}, \cite[Theorem 3.2, Theorem
4.3]{ElKaoutit/Gomez:2003}, \cite[28.21]{Brzezinski/Wisbauer:2003}
o, alternativamente, que $A$ es, como $\coring{C}$--com{\'o}dulo, un
generador simple de la categor{\'\i}a de $\coring{C}$--com{\'o}dulos por la
derecha. Vemos, por tanto, que lo que hay detr{\'a}s de \cite[2.1
Fundamental Theorem]{Sweedler:1975} puede expresarse en t{\'e}rminos
categ{\'o}ricos. De hecho, esta idea se ha explotado recientemente
para establecer una generalizaci{\'o}n de la teor{\'\i}a de Sweedler para
anillos simples artinianos \cite{Cuadra/Gomez:2005arXiv}. En este
trabajo mostramos que la idea de obtener un isomorfismo de
coanillos a partir de propiedades categ{\'o}ricas puede formularse en
{\'u}ltimo t{\'e}rmino mediante com{\'o}nadas.

Cada elemento <<como de grupo>> $g$ de un coanillo $\coring{C}$
sobre un anillo unitario $A$, da lugar \cite{Brzezinski:2002} a un
homomorfismo can\'{o}nico de $A$--coanillos $\mathsf{can}_A : A
\tensor{B} A \rightarrow \coring{C}$, que lleva $1 \tensor{B} 1$
en $g$, donde $B$ es el subanillo de los elementos
$g$--coinvariantes de $A$, y $A\tensor{B} A$ es el coanillo
can\'{o}nico de Sweedler \cite{Sweedler:1975}. El elemento <<como de
grupo>> $g$ proporciona tambi\'{e}n un par de funtores adjuntos entre
la categor\'{\i}a $\rmod{B}$ de los $B$--m\'{o}dulos por la derecha y la
categor\'{\i}a $\rcomod{\coring{C}}$ de $\coring{C}$--com\'{o}dulos por la
derecha \cite{Brzezinski:2002}. El adjunto por la izquierda est\'{a}
definido aqu\'{\i} como un producto tensor $- \tensor{B} A : \rmod{B}
\rightarrow \rcomod{\coring{C}}$, usando la estructura de
$\coring{C}$--com\'{o}dulo por la derecha que define $g$ sobre $A$. El
adjunto por la derecha puede ser interpretado como el funtor
$\hom{\coring{C}}{A}{-} : \rcomod{\coring{C}} \rightarrow
\rmod{B}$. T. Brzezi\'nski demuestra \cite[Theorem
5.6]{Brzezinski:2002} que, para $\coring{C}$ plano como
$A$--m\'{o}dulo por la izquierda, esta adjunci\'{o}n es una equivalencia
de categor\'{\i}as si, y s\'{o}lo si, $\coring{C}$ es de Galois (i.e.,
$\mathsf{can}_A$ es un isomorfismo) y $A$ es un $B$--m\'{o}dulo por la
izquierda fielmente plano. La aplicaci\'{o}n can\'{o}nica $\mathsf{can}_A$
puede ser interpretada como un homomorfismo de com\'{o}nadas de la
siguiente manera: El $A$--coanillo $\coring{C}$ da lugar a una
com\'{o}nada sobre $\rmod{A}$ construida sobre el funtor $- \tensor{A}
\coring{C}$ \cite[18.28]{Brzezinski/Wisbauer:2003}. Por otra
parte, la adjunci\'{o}n asociada a la extensi\'{o}n de anillos $B
\subseteq A$ determina \cite[Section 3.1]{Barr/Wells:2002} otra
com\'{o}nada sobre $\rmod{A}$. De esta manera, la aplicaci\'{o}n can\'{o}nica
$\mathsf{can}_A$ da lugar al homomorfismo de com\'{o}nadas $-
\tensor{A} \mathsf{can}_A : - \tensor{A} A \tensor{B} A
\rightarrow - \tensor{A} \coring{C}$. Si $(\coring{C},g)$ es un
coanillo de Galois, estas com\'{o}nadas son isomorfas y ocurre que el
funtor $- \tensor{B} A : \rmod{B} \rightarrow \rcomod{\coring{C}}$
es, salvo isomorfismos naturales, el funtor de comparaci\'{o}n de
Eilenberg-Moore \cite[Section 3.2]{Barr/Wells:2002}. Por tanto,
una de las implicaciones de \cite[Theorem 5.6]{Brzezinski:2002}
puede obtenerse como una consecuencia del Teorema de Beck
\cite[Section 3.3]{Barr/Wells:2002}. Este no parece ser el caso de
la implicaci\'{o}n rec\'{\i}proca. Concretamente, el hecho de que el funtor
$- \tensor{B} A : \rmod{B} \rightarrow \rcomod{\coring{C}}$ sea
una equivalencia implique que la aplicaci\'{o}n can\'{o}nica sea un
isomorfismo exige un argumento independiente, que reposa sobre
cierta relaci\'{o}n existente entre la aplicaci\'{o}n can\'{o}nica
$\mathsf{can}_A$ y la counidad de la adjunci\'{o}n entre los funtores
$ - \tensor{B} A$ y $\hom{\coring{C}}{A}{-}$. De hecho, para que
$\coring{C}$ sea de Galois, basta con que la counidad de dicha
adjunci\'{o}n sea un isomorfismo
\cite[18.26]{Brzezinski/Wisbauer:2003}, \cite[Lemma
3.1]{ElKaoutit/Gomez:2003}, apareciendo la condici\'{o}n de ser Galois
como parte de una caracterizaci\'{o}n del car\'{a}cter fiel y pleno del
funtor $\hom{\coring{C}}{A}{-}$
\cite[18.27]{Brzezinski/Wisbauer:2003}, \cite[Theorem
3.8]{Caenepeel/DeGroot/Vercruysse:unp2005}, \cite[Remark
3.7]{ElKaoutit/Gomez:2003}. De hecho, los citados resultados est\'{a}n
demostrados en el \'{a}mbito m\'{a}s general de los coanillos de
comatrices, introducidos en \cite{ElKaoutit/Gomez:2003}, y en los
que el papel del elemento <<como de grupo>> $g$ lo juega un
$\coring{C}$--com\'{o}dulo por la derecha $\Sigma$ que es finitamente
generado y proyectivo como $A$--m\'{o}dulo, y $B =
\rend{\coring{C}}{\Sigma}$. La generalizaci\'{o}n de \cite[Theorem
5.6]{Brzezinski:2002} en este \'{a}mbito fue demostrada en
\cite[Theorem 3.2]{ElKaoutit/Gomez:2003}. Nuestro objetivo en este
art\'{\i}culo es investigar qu\'{e} aspectos de estos resultados admiten
una formulaci\'{o}n en t\'{e}rminos puramente de com\'{o}nadas, lo que
esperamos nos permita en futuras situaciones concretas concentrar
nuestra atenci\'{o}n en lo espec\'{\i}fico de cada una de ellas, al tener
resueltos aspectos relevantes generales.

Partiremos de una com\'{o}nada $G$ sobre una categor\'{\i}a
$\cat{A}$, y de un funtor $L : \cat{B} \rightarrow \cat{A}$ con un
adjunto por la derecha $R : \cat{A} \rightarrow \cat{B}$.
Parametrizaremos de manera biun\'{\i}voca los funtores $K :
\cat{B} \rightarrow \coalg{\cat{A}}{G}$ que se factorizan a
trav\'{e}s de $L$ mediante los homomorfismos de com\'{o}nadas
$\varphi : LR \rightarrow G$ (Teorema \ref{Kphi}). Usaremos la
notaci\'{o}n $K_{\varphi}$ para designar esta dependencia.
Seguidamente, veremos bajo qu\'{e} condiciones uno de tales
funtores $K_{\varphi} : \cat{B} \rightarrow \coalg{\cat{A}}{G}$
tiene un adjunto por la derecha $D_{\varphi} : \coalg{\cat{A}}{G}
\rightarrow \cat{B}$ (Proposici\'{o}n \ref{adjuncion}).
Demostraremos que $D_{\varphi}$ es fiel y pleno si y s\'{o}lo si
$\varphi$ es un isomorfismo y $L$ preserva ciertos igualadores
(Teorema \ref{Descent}) y concluiremos nuestros resultados
generales caracterizando cu\'{a}ndo $K_{\varphi}$ proporciona una
equivalencia entre las categor\'{\i}as $\cat{B}$ y
$\coalg{\cat{A}}{G}$ (Teorema \ref{comonadic}). Obviamente, los
funtores caracterizados as\'{\i} son, a fortiori, comon\'{a}dicos
o tripeables pero, a diferencia del planteamiento del Teorema de
Beck, aqu\'{\i} la com\'{o}nada $G$ est\'{a} dada de antemano, y
cada funtor $K_{\varphi}$ corresponde a una <<representaci\'{o}n>>
de $G$. La situaci\'{o}n tratada por el Teorema de Beck es aquella
en que $\varphi$ es la identidad, esto es, $G = LR$ es la
com\'{o}nada asociada a la adjunci\'{o}n.

Una situaci\'{o}n que motiva nuestro punto de vista viene dada por una
estructura entrelazante entre un \'{a}lgebra y una co\'{a}lgebra, y un
m\'{o}dulo entrelazado \cite{Brzezinski/Majid:1998}. Entonces la
com\'{o}nada $G$ viene dada por el coanillo asociado a la estructura
entrelazante \cite[Section 2]{Brzezinski:2002}, y el funtor
$K_{\varphi}$ est\'{a} definido por un m\'{o}dulo entralazado, que no es
sino un com\'{o}dulo sobre el mencionado coanillo. Una cuesti\'{o}n natural
es estudiar la relaci\'{o}n entre la categor\'{\i}a de m\'{o}dulos entrelazados y
la de m\'{o}dulos sobre el subanillo de coinvariantes del m\'{o}dulo
entrelazado y, en particular, si ambas categor\'{\i}as son equivalentes.

Aplicaremos nuestros teoremas generales al caso de coanillos sobre
anillos firmes, lo que ilustrar\'{a} c{\'o}mo los resultados sobre
com\'{o}nadas de las secciones \ref{candefinido} y \ref{adjeq}
simplifican significativamente el tratamiento de algunos aspectos
relevantes de los coanillos de comatrices y los com\'{o}dulos de
Galois estudiados en \cite{ElKaoutit/Gomez:2003},
\cite{Gomez/Vercruysse:2005arXiv}.

\medskip

\noindent \textbf{Nota importante:} En esta versi{\'o}n revisada,
debemos decir que buena parte de los resultados de las secciones
\ref{candefinido} y \ref{adjeq} son conocidos para los
especialistas en Teor{\'\i}a de Categor{\'\i}as. De hecho, B. Mesablishvili
nos ha informado gentilmente de este hecho, proporcion{\'a}ndonos en
particular las referencias \cite{Dubuc:1970} y \cite{Beck:2003}.
Hemos incluido las correspondientes atribuciones de los
resultados, aunque no se puede descartar que m{\'a}s enunciados de las
secciones \ref{candefinido} y \ref{adjeq} trabajo resulten, al
menos, familiares para los mencionados especialistas. Esperamos,
no obstante, que nuestros enunciados y pruebas elementales
pudieran ser de alguna utilidad a aquellos lectores con
conocimientos no especializados en Categor{\'\i}as y Funtores, al hacer
m{\'a}s accesible para ellos, como as{\'\i} ha sido para el autor de estas
notas, un enlace entre los coanillos de Galois y la bien
desarrollada teor{\'\i}a general de com{\'o}nadas.

\section{Funtores con valores en co\'{a}lgebras y homomorfismos de
com\'{o}nadas}\label{candefinido}

Sea $(G,\Delta, \varepsilon)$ una com\'{o}nada (o cotriple) sobre una
categor\'{\i}a $\cat{A}$, esto es, un funtor $G : \cat{A} \rightarrow
\cat{A}$ junto con dos transformaciones naturales $\Delta : G
\rightarrow G^2$ y $\varepsilon : G \rightarrow id_{\cat{A}}$
tales que los diagramas
\begin{equation*}
\xymatrix{G \ar^{\Delta}[r] \ar_{\Delta}[d] & G^2 \ar^{G\Delta}[d]
\\ G^2 \ar_{\Delta G}[r] & G^3 } \qquad \xymatrix{G & G^2 \ar_{\varepsilon G}[l] \ar^{G \varepsilon}[r] & G \\
& \ar@{=}[ul] G \ar|{\Delta}[u] \ar@{=}[ur] & }
\end{equation*}
son conmutativos \cite[Chapter 3]{Barr/Wells:2002}. Seguiremos en
lo posible \cite{Barr/Wells:2002}, entendiendo autom\'{a}ticamente
cada afirmaci\'{o}n sobre m\'{o}nadas (triples) en su versi\'{o}n para
com\'{o}nadas. Supongamos un funtor $L : \cat{B} \rightarrow \cat{A}$
que tiene un adjunto por la derecha $R : \cat{A} \rightarrow
\cat{B}$. Entonces, si $\eta : id_{\cat{B}} \rightarrow RL$ es la
unidad de la adjunci\'{o}n, y $\epsilon : LR \rightarrow id_{\cat{A}}$
es su counidad, se tiene definida \cite[Proposition
3.1.2]{Barr/Wells:2002} la com\'{o}nada $(LR, \delta, \epsilon)$ sobre
$\cat{A}$, donde $\delta = L\eta R$. Recordemos la categor\'{\i}a
$\coalg{\cat{A}}{G}$ de las $G$--co\'{a}lgebras \cite[Section
3.1]{Barr/Wells:2002}, cuyos objetos son pares $(X,x)$
consistentes en un objeto $X$ de $\cat{A}$ y un morfismo $x : X
\rightarrow GX$ tales que
\begin{equation*}
Gx \circ x = \Delta_X \circ x, \qquad \varepsilon_X \circ x = id_X
\end{equation*}
Dadas $G$--co\'{a}lgebras $(X,x), (X',x')$, los homomorfismos $f : X
\rightarrow X'$ en $\cat{A}$ tales que $Gf \circ x = x' \circ f$
constituyen el conjunto de homomorfismos
$\hom{\coalg{\cat{A}}{G}}{X}{X'}$ en $\coalg{\cat{A}}{G}$ de
$(X,x)$ a $(X',x')$.

Denotemos por $U : \coalg{\cat{A}}{G} \rightarrow \cat{A}$ el
funtor que olvida. Consideraremos aquellos funtores $K : \cat{B}
\rightarrow \coalg{\cat{A}}{G}$ tales que el diagrama
\begin{equation}\label{UKL}
\xymatrix{\cat{B} \ar^{K}[r] \ar_{L}[dr] & \coalg{\cat{A}}{G}
\ar^{U}[d] \\
& \cat{A}}
\end{equation}
es conmutativo. Estamos interesados particularmente en el caso en
que $K$ proporciona una equivalencia de categor\'{\i}as entre
$\cat{B}$ y $\coalg{\cat{A}}{G}$. Comenzaremos estableciendo una
correspondencia biun\'{\i}voca entre los funtores $K$ que hacen
conmutar el diagrama \eqref{UKL} y los homomorfismos de
com\'{o}nadas $\varphi : LR \rightarrow G$. Esta correspondencia
es consecuencia de la siguiente Proposici\'{o}n
\ref{correspondencia} que, por otra parte, establece algunos
hechos t\'{e}cnicos fundamentales para abordar la
caracterizaci\'{o}n de las equivalencias de categor\'{\i}as $K$
que hacen conmutar \eqref{UKL}. Recordemos \cite[Section
3.6]{Barr/Wells:2002} que un homomorfismo de com\'{o}nadas de $LR$
a $G$ es una transformaci\'{o}n natural $\varphi : LR \rightarrow
G$ tal que $\Delta \varphi = \varphi^2 \delta$ y $\varepsilon
\varphi = \epsilon$. La correspondencia biun{\'\i}voca entre los
objetos descritos en las afirmaciones (A) y (C) de la Proposici{\'o}n
\ref{correspondencia} puede deducirse de \cite[Proposition
II.1.4]{Dubuc:1970}.

\begin{proposition}\label{correspondencia}
Existe una correspondencia biun\'{\i}voca entre \\

\noindent (A) Morfismos de com\'{o}nadas de $(LR, L\eta R, \epsilon)$
a
$(G,\Delta,\varepsilon)$, \\

\medskip

 \noindent (B) Transformaciones naturales $\xymatrix{R
\ar^{\alpha}[r] & RG}$ tales que los siguientes diagramas conmutan
\begin{equation}\label{comod}
\xymatrix{R \ar^{\alpha}[r] \ar^{\alpha}[d] & RG \ar^{R\Delta}[d] \\
RG \ar^{\alpha G}[r] & RG^2 } \qquad \xymatrix{ R \ar^{\alpha}[r]
\ar@{=}[dr] & RG \ar^{R \varepsilon}[d] \\
& R }
\end{equation}
y

\noindent (C) Transformaciones naturales $\xymatrix{L
\ar^{\beta}[r] & GL}$ tales que los siguientes diagramas conmutan
\begin{equation}\label{comodb}
\xymatrix{L \ar^{\beta}[r] \ar^{\beta}[d] & GL \ar^{\Delta L}[d] \\
GL \ar^{ G \beta }[r] & G^2L } \qquad \xymatrix{ L \ar^{\beta}[r]
\ar@{=}[dr] & GL \ar^{\varepsilon L}[d] \\
& L }
\end{equation}
\end{proposition}
\begin{proof}
Comenzamos demostrando la correspondencia biun\'{\i}voca entre las
transformaciones naturales descritas en $(A)$ y $(B)$.  Si
partimos de un morfismo de com\'{o}nadas $\varphi : LR \rightarrow G$,
entonces podemos definir la transformaci\'{o}n natural
\begin{equation}\label{rhophi}
\xymatrix{ R \ar@/_1pc/_{\alpha}[rr] \ar^{\eta R}[r] & RLR \ar^{R
\varphi}[r] & RG}
\end{equation}
Para comprobar que conmuta el primer diagrama en \eqref{comod},
consideremos un objeto $X$ de $\cat{A}$ y calculamos
\begin{equation*}
\begin{array}{lclr}
R\Delta_X \circ \alpha_X & = & R\Delta_X \circ R\varphi_X \circ
\eta_{RX} & \\
& = & R\varphi^2_X \circ RL\eta_{RX} \circ \eta_{RX} & (\varphi
\hbox{ es de com\'{o}nadas}) \\
 & = & R\varphi_{GX} \circ RLR\varphi_X \circ RL\eta_{RX} \circ
 \eta_{RX} & (\varphi^2_X = \varphi_{GX} \circ LR\varphi_X) \\
 & = & R\varphi_{GX} \circ RLR\varphi_X \circ \eta_{RLRX} \circ
 \eta_{RX} & (\eta \hbox{ es natural}) \\
 & = & R\varphi_{GX} \circ \eta_{RGX} \circ R\varphi_X \circ \eta_{RX} &
 (\eta \hbox{ es natural}) \\
 & = & \alpha_{GX} \circ \alpha_X
\end{array}
\end{equation*}
Para el segundo diagrama en \eqref{comod}, tenemos
\begin{equation*}
\begin{array}{lclr}
R\varepsilon_X \circ \alpha_X & = & R\varepsilon_X \circ
R\varphi_X
\circ \eta_{RX} & \\
& = & R\epsilon_X \circ \eta_{RX} & (\varphi \hbox{ es de
com\'{o}nadas}) \\
& = & id_{RX} & (\hbox{por adjunci\'{o}n})
\end{array}
\end{equation*}
En sentido inverso, partiendo de una transformaci\'{o}n natural
$\alpha : R \rightarrow RG$ que verifique \eqref{comod}, definimos
la transformaci\'{o}n natural
\begin{equation}\label{can}
\xymatrix{LR \ar^{L\alpha}[r] \ar@/_1pc/_{\varphi}[rr] & LRG
\ar^{\epsilon G}[r] & G}
\end{equation}
Para demostrar que $\varphi$, definida en \eqref{can}, es un
homomorfismo de com\'{o}nadas, necesitaremos utilizar que, para cada
objeto $X$ de $\cat{A}$, se tienen las igualdades
\begin{equation}\label{can2}
\varphi^2_X = \epsilon_{G^2X} \circ L\alpha_{GX} \circ
LR\epsilon_{GX} \circ LRL\alpha_X
\end{equation}
y
\begin{equation*}
\varphi^2_X = G\epsilon_{GX} \circ GL\alpha_X \circ
\epsilon_{GLRX} \circ L\alpha_{RLX},
\end{equation*}
por definici\'{o}n de $\varphi$. Realizamos el siguiente c\'{a}lculo:
\begin{equation*}
\begin{array}{lclr}
\varphi^2_X \circ L\eta_{RX} & = & \epsilon_{G^2X} \circ
L\alpha_{GX} \circ LR\epsilon_{GX} \circ LRL\alpha_X \circ
L\eta_{RX}
& (\hbox{por \eqref{can2}}) \\
& = &  \epsilon_{G^2X} \circ L\alpha_{GX} \circ LR\epsilon_{GX}
\circ L\eta_{RGX} \circ L\alpha_X & (\eta \hbox{ es natural}) \\
& = & \epsilon_{G^2X} \circ L\alpha_{GX}  \circ L\alpha_X &
(R\epsilon_{GX} \circ \eta_{RGX} = id_{RGX}) \\
& = & \epsilon_{G^2X} \circ LR\Delta_X \circ L\alpha_X &
(\alpha_{GX} \circ \alpha_X = R\Delta_X \circ \alpha_X) \\
& = & \Delta_X \circ \epsilon_{GX} \circ L\alpha_X & (\epsilon
\hbox{ es natural}) \\
& = & \Delta_X \circ \varphi_X
\end{array}
\end{equation*}
Para comprobar la segunda condici\'{o}n que define un homomorfismo de
com\'{o}nadas, tenemos
\begin{equation*}
\begin{array}{lclr}
\varepsilon_X \circ \varphi_X & = & \varepsilon_X \circ
\epsilon_{GX} \circ L\alpha_X & \\
& = & \epsilon_X \circ LR\varepsilon_X \circ L\alpha_X & (\epsilon
\hbox{ es natural}) \\
& = & \epsilon_X & (R\varepsilon_X \circ \alpha_X = id_{RX})
\end{array}
\end{equation*}
Bien, partamos ahora de un homomorfismo de com\'{o}nadas $\varphi
: LR \rightarrow G$ y construyamos la transformaci\'{o}n natural
$\alpha : R \rightarrow RG$ de acuerdo con \eqref{rhophi}.
Consideremos ahora el homomorfismo de com\'{o}nadas, digamos
$\varphi'$, definido, a partir de tal $\alpha$, en \eqref{can}, y
veamos que coincide con el morfismo original $\varphi$. Para ello,
calculemos, para un objeto $X$ de $\cat{A}$,
\begin{equation*}
\begin{array}{lclr}
\varphi'_X & = & \epsilon_{GX} \circ L\alpha_X & \\
& = & \epsilon_{GX} \circ LR\varphi_X \circ L\eta_{RX} & \\
& = & \varphi_X \circ \epsilon_{LRX} \circ L\eta_{RX} & (\epsilon
\hbox{ es natural}) \\
& = & \varphi_X &
\end{array}
\end{equation*}
Por \'{u}ltimo hemos de ver que, partiendo de una
transformaci\'{o}n natural $\alpha : R \rightarrow RG$ que
satisfaga \eqref{comod} y construyendo sucesivamente el morfismo
de com\'{o}nadas $\varphi$ seg\'{u}n \eqref{can}, y la nueva
transformaci\'{o}n natural, digamos $\alpha'$, de acuerdo con
\eqref{rhophi}, volvemos a obtener el $\alpha$ original. Esto se
sigue del siguiente c\'{a}lculo, para $X$ un objeto de $\cat{A}$:
\begin{equation*}
\begin{array}{lclr}
\alpha'_X & = & R\varphi_X \circ \eta_{RX} & \\
& = & R\epsilon_{GX} \circ RL\alpha_X \circ \eta_{RX} & (\eta
\hbox{
es natural}) \\
& = & R\epsilon_{GX} \circ \eta_{RGX} \circ \alpha_X & \\
& = & \alpha_X &
\end{array}
\end{equation*}
Abordemos ahora la demostraci\'{o}n de que hay una correspondencia
biun\'{\i}voca entre las transformaciones descritas en $(A)$ y las
descritas en $(C)$. As\'{\i}, si partimos de un homomorfismo de
com\'{o}nadas $\varphi: LR \rightarrow G$, definimos la transformaci\'{o}n
natural
\begin{equation}\label{lambdaphi}
\xymatrix{ L \ar@/_1pc/_{\beta}[rr] \ar^{L \eta}[r] & LRL \ar^{
\varphi L}[r] & GL}
\end{equation}
y comprobemos que verifica que los diagramas \eqref{comodb}
conmutan. Tomamos, pues, un objeto $Y$ de $\cat{B}$ y calculamos
\begin{equation*}
\begin{array}{lclr}
\Delta_{LY} \circ \beta_Y & = & \Delta_{LY} \circ \varphi_{LY}
\circ L\eta_Y & \\
& = & \varphi_{LY}^2 \circ L\eta_{RLY} \circ L\eta_Y & (\varphi
\hbox{ es de com\'{o}nadas}) \\
& = & G\varphi_{LX} \circ \varphi_{LRLY} \circ L\eta_{RLY} \circ
L\eta_Y & (\varphi^2_{LY} = G\varphi_{LX} \circ \varphi_{LRLY}) \\
& = & G\varphi_{LY} \circ \varphi_{RLRY} \circ LRL\eta_Y \circ
L\eta_Y & (\eta \hbox{ es natural}) \\
& = & G\varphi_{LY} \circ GL\eta_Y \circ \varphi_{LY} \circ
L\eta_Y & (\varphi \hbox{ es natural}) \\
& = & G\beta_Y \circ \beta_Y
\end{array}
\end{equation*}
Para comprobar que es segundo diagrama de \eqref{comodb} conmuta,
hacemos el siguiente c\'{a}lculo:
\begin{equation*}
\begin{array}{lclr}
\varepsilon_{LY} \circ \beta_Y & = & \varepsilon_{LY} \varphi_{LY}
\circ L\eta_Y & \\
& = & \epsilon_{LY} \circ L\eta_Y & (\varphi \hbox{ es de
com\'{o}nadas}\\
& = & id_{LY} & (\hbox{por adjunci\'{o}n})
\end{array}
\end{equation*}
Rec\'{\i}procamente, a cada transformaci\'{o}n natural $\beta : L
\rightarrow GL$ que haga conmutar los diagramas \eqref{comodb}, le
asignamos la transformaci\'{o}n natural $\varphi$ definida por
\begin{equation}\label{canr}
\xymatrix{LR \ar^{\beta R}[r] \ar@/_1pc/_{\varphi}[rr] & GLR \ar^{
G \epsilon}[r] & G}
\end{equation}
que vamos a comprobar seguidamente es un homomorfismo de
com\'{o}nadas. Para la primera condici\'{o}n requerida para ser un
homomorfismo de com\'{o}nadas calculamos, para un objeto $X$ de
$\cat{A}$:
\begin{equation*}
\begin{array}{lclr}
\varphi^2_X \circ L\eta_{RX} & = & G\varphi_X \circ \varphi_{LRX}
\circ L\eta_{RX} & \\
& = & G^2\epsilon_X \circ G \beta_{RX} \circ G \epsilon_{LRX}
\circ \beta_{RLRX} \circ L\eta_{RX} & \\
& = & G^2\epsilon_X \circ G\beta_{RX} \circ G\epsilon_{LRX} \circ
GL\eta_{RX} \circ \beta_{RX} & (\beta \hbox{ es natural}) \\
& = & G^2\epsilon_{X} \circ G \beta_{RX} \circ \beta_{RX} &  (\epsilon_{LRX} \circ L\eta_{RX} = id_{LRX})\\
& = & G^2\epsilon_X \circ \Delta_{LRX} \circ \beta_{RX} &
(\hbox{por \eqref{comodb}}) \\
& = & \Delta_X \circ G\epsilon_X \circ \beta_{RX} & (\Delta \hbox{
es natural}) \\
& = & \Delta_X \circ \varphi_X &
\end{array}
\end{equation*}
Una comprobaci\'{o}n de la segunda condici\'{o}n para ser homomorfismo de
com\'{o}nadas es la siguiente:
\begin{equation*}
\begin{array}{lclr}
\varepsilon_X \circ \varphi_X & = & \varepsilon_X \circ G
\epsilon_X \circ \beta_{RX} & \\
& = & \epsilon_X \circ \varepsilon_{LRX} \circ \beta_{RX} &
(\varepsilon \hbox{ es natural}) \\
& = & \epsilon_X & (\hbox{por \eqref{comodb}})
\end{array}
\end{equation*}
Veamos por \'{u}ltimo que estas correspondencias son mutuamente
inversas. As\'{\i} pues, partiendo de un homomorfismo de com\'{o}nadas
$\varphi : LR \rightarrow G$, y siendo $\beta$ la transformaci\'{o}n
natural definida en \eqref{lambdaphi}, tenemos, para $X$ un objeto
de $\cat{A}$, y $\varphi'$ el homomorfismo de com\'{o}nadas definido a
partir de $\beta$ seg\'{u}n \eqref{canr}:
\begin{equation*}
\begin{array}{lclr}
\varphi'_X & = & G \epsilon_X \circ \beta_{RX} & \\
& = & G\epsilon_X \circ \varphi_{LRX} \circ L\eta_{RX} &
\\
& = & \varphi_X \circ LR\epsilon_X \circ L\eta_{RX} & (\varphi
\hbox{ es natural}) \\
& = & \varphi_X & (R \epsilon_X \circ \eta_{RX} = id_{RX})
\end{array}
\end{equation*}
Y partiendo de una transformaci\'{o}n natural $\beta : L
\rightarrow GL$ sujeta a las condiciones \eqref{comodb}, tenemos
definido un homomorfismo de com\'{o}nadas $\varphi$ seg\'{u}n
\eqref{canr}, que permite definir, de acuerdo con
\eqref{lambdaphi}, una transformaci\'{o}n natural $\beta'$.
Tenemos, para cada objeto $Y$ de $\cat{B}$:
\begin{equation*}
\begin{array}{lclr}
\beta'_Y & = & \varphi_{LY} \circ L\eta_Y & \\
& = & G \epsilon_{LY} \circ \beta_{RLY} \circ L\eta_Y & \\
& = & \beta_Y \circ \epsilon_{LY} \circ L\eta_Y & (\beta \hbox{ es
natural}) \\
& = & \beta_Y
\end{array}
\end{equation*}
Esto concluye la prueba.
\end{proof}

Conviene que nos dotemos de una terminolog\'{\i}a c\'{o}moda para
referirnos a la situaci\'{o}n descrita por la Proposici\'{o}n
\ref{correspondencia}. As\'{\i}, si el dato de partida es un
homomorfismo de com\'{o}nadas $\varphi : LR \rightarrow G$,
entonces nos referiremos a las transformaciones naturales $\alpha
: R \rightarrow RG$ y $\beta : L \rightarrow GL$ como la
\emph{representaci\'{o}n co-inducida} (resp.
\emph{representaci\'{o}n inducida}) de $\varphi$. En su conjunto,
hablaremos de las \emph{representaciones} de $\varphi$. Cuando
partamos de una transformaci\'{o}n natural $\alpha$ sujeta a las
condiciones \eqref{comod}, o de una transformaci\'{o}n natural
$\beta$ sujeta a \eqref{comodb}, entonces el homomorfismo de
com\'{o}nadas correspondiente $\varphi : LR \rightarrow G$ se
llamar\'{a} \emph{transformaci\'{o}n can\'{o}nica} asociada a
$\alpha$ (resp. $\beta$).

Dado un homomorfismo de com\'{o}nadas $\varphi : LR \rightarrow G$
es \'{u}til escribir las ecuaciones que lo ligan con sus
representaciones asociadas $\alpha : R \rightarrow RG$ y $\beta :
L \rightarrow GL$. Concretamente, de la demostraci\'{o}n de la
Proposici\'{o}n \ref{correspondencia}, deducimos que
\begin{equation}\label{tabladef}
\begin{array}{ccc}
\alpha = R \varphi \circ \eta R & & \varphi = \epsilon G \circ
L\alpha \\
 & \\
\beta = \varphi L \circ L \eta & & \varphi = G \epsilon \circ
\beta R
\end{array}
\end{equation}

Pasamos a deducir de la Proposici\'{o}n \ref{correspondencia} el
resultado fundamental de esta secci\'{o}n, que puede tambi{\'e}n
obtenerse de \cite[Theorem II.1.1]{Dubuc:1970}. Por otra parte, B.
Mesablishvili nos ha comunicado que el Teorema \ref{Kphi} puede
deducirse de cierta meta-adjunci{\'o}n, descrita en \cite[p{\'a}g.
11]{Beck:2003}, y para cuya prueba remite el mismo Beck a
\cite{Eilenberg/Moore:1965}, entre ciertas meta-categor{\'\i}as
adecuadas, de una parte $Adj(\cat{A})$, cuyos objetos son
adjunciones, y de otra $Trip(\cat{A})^{op}$, cuyos objetos son
m{\'o}nadas (o triples), siempre con los morfismos adecuados.

\begin{theorem}\label{Kphi}
Dada una com\'{o}nada $G$ sobre $\cat{A}$ y un funtor $L : \cat{B}
\rightarrow \cat{A}$, si $L$ tiene un adjunto por la derecha $R :
\cat{A} \rightarrow \cat{B}$, entonces existe una correspondencia
biyectiva entre funtores $K : \cat{B} \rightarrow
\coalg{\cat{A}}{G}$ que hacen conmutar \eqref{UKL} y homomorfismos
de com\'{o}nadas $\varphi : LR \rightarrow G$.
\end{theorem}
\begin{proof}
Dado que, seg\'{u}n la Proposici\'{o}n \ref{correspondencia}, existe una
correspondencia biun\'{\i}voca entre homomorfismos de com\'{o}nadas
$\varphi : LR \rightarrow G$ y transformaciones naturales $\beta :
L \rightarrow GL$ que hagan conmutar \eqref{comodb}, basta con que
demostremos que estas \'{u}ltimas est\'{a}n en correspondencia biun\'{\i}voca
con los funtores $K : \cat{B} \rightarrow \coalg{\cat{A}}{G}$
tales que $UK = L$. Esta correspondencia va como sigue: dada
$\beta : L \rightarrow GL$, definimos $KY = (LY, \beta_Y)$ para
cada objeto $Y$ de $\cat{B}$, y $Kf = Lf$ para cada morfismo $f$
en $\cat{B}$. Los diagramas \eqref{comodb} muestran entonces que
$(LY,\beta_Y)$ es una $G$--co\'{a}lgebra y la naturalidad de $\beta$
implica que $Lf$ es un homomorfismo de $G$--co\'{a}lgebras.
Rec\'{\i}procamente, dado el funtor $K : \cat{B} \rightarrow
\coalg{\cat{A}}{G}$ definimos, para cada objeto $Y$ de $\cat{B}$,
$\beta_Y : LY \rightarrow GLY$ como el homomorfismo estructura de
la $G$--co\'{a}lgebra $KY$. Es f\'{a}cil comprobar que esta asignaci\'{o}n
define una transformaci\'{o}n natural $\beta : L \rightarrow GL$ que
hace conmutar \eqref{comodb}.
\end{proof}

El Teorema \ref{Kphi} puede ser considerado como una
generalizaci\'{o}n de \cite[Theorem 3.3]{Barr/Wells:2002}.

\begin{corollary}\cite[Theorem 3.3]{Barr/Wells:2002}
Sean $H$ y $G$ com\'{o}nadas sobre una categor\'{\i}a $\cat{A}$. Existe una
correspondencia biun\'{\i}voca entre homomorfismos de com\'{o}nadas
$\varphi : H \rightarrow G$ y funtores $K : \coalg{\cat{A}}{H}
\rightarrow \coalg{\cat{A}}{G}$ que hacen conmutar el diagrama
\begin{equation*}
\xymatrix{\coalg{\cat{A}}{H} \ar^{K}[rr] \ar_{U_H}[dr]& &
\coalg{\cat{A}}{G} \ar^{U_G}[dl]
\\ & \cat{A} & }
\end{equation*}
donde $U_H$ y $U_G$ son los funtores que olvidan.
\end{corollary}
\begin{proof}
El funtor $U_H$ tiene un adjunto por la derecha $F_H$ tal que la
com\'{o}nada asociada es $H$ \cite[Secci\'{o}n 3.2]{Barr/Wells:2002}. El
corolario se sigue ahora del Teorema \ref{Kphi} tomando $\cat{B} =
\coalg{\cat{A}}{H}$ y $L = U_H$.
\end{proof}

\section{Adjunciones y equivalencias}\label{adjeq}

Partimos de una com\'{o}nada $(G,\Delta,\varepsilon)$ sobre una
categor\'{\i}a $\cat{A}$, y $L : \cat{B} \rightarrow \cat{A}$ un
funtor. Supondremos, adem\'{a}s, que $L$ tiene un adjunto por la
derecha $R : \cat{A} \rightarrow \cat{B}$. De acuerdo con el
Teorema \ref{Kphi}, los funtores $K : \cat{B} \rightarrow
\coalg{\cat{A}}{G}$ tales que $UK = L$, para $U :
\coalg{\cat{A}}{G} \rightarrow \cat{A}$ el funtor que olvida,
est\'{a}n en correspondencia biun\'{\i}voca con los homomorfismos de
com\'{o}nadas $\varphi : LR \rightarrow G$. Haremos patente esta
dependencia escribiendo $K_{\varphi}$ para el funtor
correspondiente a $\varphi$. Expl\'{\i}citamente, el funtor
$K_{\varphi}$ est\'{a} definido como
\begin{equation*}
\xymatrix{K_\varphi : \cat{B} \ar[r] & \coalg{\cat{A}}{G} & & (Y
\mapsto (LY, \varphi_{LY} \circ L\eta_Y)),}
\end{equation*}
sobre objetos, y como $L$ sobre homomorfismos. Tendremos presentes
la Proposici\'{o}n \ref{correspondencia} y las ecuaciones
\eqref{tabladef} que relacionan $\varphi$ con sus representaciones
$\alpha : R \rightarrow RG$ y $\beta : L \rightarrow GL$.

\medskip

Comenzamos estudiando cu\'{a}ndo $K_{\varphi}$ tiene un adjunto
por la derecha. La siguiente proposici\'{o}n es un caso particular
de \cite[Theorem A.1]{Dubuc:1970}. Damos una prueba elemental en
nuestro caso.

\begin{proposition}\label{adjuncion}
Supongamos que para cada $G$--co\'{a}lgebra $(X,x)$ existe en
$\cat{B}$ el igualador del par de morfismos $\alpha_X, Rx : RX
\rightarrow RGX$. Entonces el funtor $K_{\varphi} : \cat{B}
\rightarrow \coalg{\cat{A}}{G}$ tiene un adjunto por la derecha
$D_{\varphi} : \coalg{\cat{A}}{G} \rightarrow \cat{B}$, cuyo valor
en $(X,x)$ es el igualador
\begin{equation}\label{Dphi}
\xymatrix{D_{\varphi}X \ar^{eq_X}[rr] & & RX
\ar@<.5ex>^-{\alpha_X}[rr] \ar@<-.5ex>_-{Rx}[rr] & & RGX}
\end{equation}
\end{proposition}
\begin{proof}
Dados objetos $X$ de $\cat{A}$ e $Y$ de $\cat{B}$, denotemos por
$\Phi : \hom{\cat{A}}{LY}{X} \rightarrow \hom{\cat{B}}{Y}{RX}$ el
isomorfismo de adjunci\'{o}n, que, en t\'{e}rminos de la counidad $\eta$,
viene definido por $\Phi (h) = Rh \circ \eta_Y$ para $h : LY
\rightarrow X$. Consideremos el siguiente diagrama conmutativo de
aplicaciones entre conjuntos:
\begin{equation}\label{adjdia}
\xymatrix{\hom{\coalg{\cat{A}}{G}}{K_{\varphi}Y}{X}
\ar@{-->}^{\Phi}[rr] \ar[d] & &
\hom{\cat{B}}{Y}{D_{\varphi}X} \ar[d] \\
\hom{\cat{A}}{LY}{X} \ar^{\Phi}[rr] \ar@<-1ex>_{x \circ - }[d] \ar@<1ex>^{G(-) \circ \beta_X}[d]& &
\hom{\cat{B}}{Y}{RX}  \ar@<-1ex>_{\hom{\cat{B}}{Y}{Rx}}[d] \ar@<1ex>^{\hom{\cat{B}}{Y}{\alpha_X}}[d]\\
\hom{\cat{A}}{LY}{GX} \ar^{\Phi}[rr]
 & &
\hom{\cat{B}}{Y}{RGX},}
\end{equation}
donde la flecha discontinua no est\'{a} a\'{u}n definida. Comprobemos que
el cuadrado inferior conmuta serialmente, en el sentido de
\cite[page 112]{Barr/Wells:2002}: partiendo de $h : LY \rightarrow
X$, tenemos
\begin{equation*}
\begin{array}{lcl}
\left[ \Phi \circ (x \circ -) \right](h) & = & \Phi(x \circ h) \\
& = & Rx \circ Rh \circ \eta_Y \\
& = & \hom{\cat{B}}{Y}{Rx}(Rh \circ \eta_Y) \\
& = & \left[\hom{\cat{B}}{Y}{Rx} \circ \Phi \right](h)
\end{array}
\end{equation*}
Por otra parte,
\begin{equation*}
\begin{array}{lclr}
\left[ \Phi \circ (G(-) \circ \beta_Y)\right](h) & = & \Phi (Gh
\circ \beta_Y) & \\
& = & RGh \circ R\beta_Y \circ \eta_Y & \\
& = & RGh \circ R\varphi_{LY} \circ RL\eta_{Y} \circ \eta_Y &
(\beta_Y = \varphi_{LY} \circ L\eta_Y) \\
& = & RGh \circ R\varphi_{LY} \circ \eta_{RLY} \circ \eta_Y &
(\eta \hbox{ es natural}) \\
& = & R\varphi_X \circ RLRh \circ \eta_{RLY} \circ \eta_Y &
(\varphi \hbox{ es natural}) \\
& = & R\varphi_X \circ \eta_{RX} \circ Rh \circ \eta_Y & (\eta
\hbox{ es natural}) \\
& = & \alpha_X \circ Rh \circ \eta_Y & \\
& = & \hom{\cat{B}}{Y}{\alpha_X}(Rh \circ \eta_Y) & \\
& = & \left[ \hom{\cat{B}}{Y}{\alpha_X} \circ \Phi \right](h).
\end{array}
\end{equation*}
Ahora bien, los dos lados verticales son igualadores, el de la
izquierda por definici\'{o}n de homomorfismo de $G$--co\'{a}lgebras, y el
de la derecha, por la propiedad universal del igualador
\eqref{Dphi}. Por tanto, el isomorfismo natural $\Phi :
\hom{\cat{A}}{LY}{X} \rightarrow \hom{\cat{B}}{Y}{RX}$ induce, por
restricci\'{o}n, un isomorfismo natural $\Phi :
\hom{\coalg{\cat{A}}{G}}{K_\varphi Y}{X} \rightarrow
\hom{\cat{B}}{Y}{D_{\varphi}X}$.
\end{proof}

\begin{remark}\label{unitcounit}
Si tomamos en el diagrama \eqref{adjdia} $X = K_{\varphi}Y$, para
un objeto $Y$ de $\cat{B}$, dado que $id_{K_{\varphi}Y}$ es un
homomorfismo de $G$--co\'{a}lgebras, deducimos que $\eta_Y =
\Phi(id_{K_{\varphi}Y})$ se factoriza a trav\'{e}s de $D_{\varphi}Y$,
de manera que la unidad de la adjunci\'{o}n $K_{\varphi} \dashv
D_{\varphi}$, denotada por $\widehat{\eta}$, viene determinada
un\'{\i}vocamente en $Y$ por la propiedad universal de un igualador,
seg\'{u}n el siguiente diagrama:
\begin{equation}\label{unitphi} \xymatrix{D_{\varphi}K_{\varphi}Y \ar[r] &
RLY \ar@<.5ex>^-{\alpha_{LY}}[rr] \ar@<-.5ex>_-{R\beta_Y}[rr] & &
RGLY \\
& Y \ar@{-->}^{\widehat{\eta}_Y}[ul] \ar^{\eta_Y}[u] & & }
\end{equation}
Para hacer expl\'{\i}cita la counidad $\widehat{\epsilon}$ en un objeto
$(X,x)$ de $\coalg{\cat{A}}{G}$, aplicamos el funtor $K_{\varphi}$
al igualador \eqref{Dphi} y obtenemos el diagrama
\begin{equation}\label{counidad1}
\xymatrix{K_{\varphi}D_{\varphi}X \ar^{Leq_X}[rr]
\ar@{-->}^{\widehat{\epsilon}_X}[drr] & & LRX
\ar@<.5ex>^-{L\alpha_X}[rr] \ar@<-.5ex>_-{LRx}[rr] \ar^{\epsilon_X}[d] & & LRGX \\
& & X, & & }
\end{equation}
es decir, $\widehat{\epsilon}_X = \epsilon_X \circ Leq_X$.
\end{remark}

Nuestro pr\'{o}ximo objetivo es demostrar que $D_{\varphi}$ es un
funtor pleno y fiel si, y s\'{o}lo si, $\varphi$ es in isomorfismo de
m\'{o}nadas y el funtor $L$ preserva los igualadores \eqref{Dphi}.
Trabajaremos primero para una $G$--co\'{a}lgebra $X$, y luego
argumentaremos globalmente. Supongamos que existe en $\cat{A}$ el
igualador
\begin{equation}\label{ecuprima}
\xymatrix{EX \ar^{eq'_X}[rr] & & LRX \ar@<.5ex>^-{L\alpha_X}[rr]
\ar@<-.5ex>_-{LRx}[rr] & & LRGX}
\end{equation}
Consideramos el siguiente diagrama en $\cat{A}$:
\begin{equation}\label{clave}
\xymatrix{ & & & &   LRGX \ar^-{\varphi_{GX}}[dd] \\
& & &   LRX \ar@<-.5ex>_{L\alpha_X}[ru] \ar@<.5ex>^{LRx}[ru] \ar^-{\varphi_{X}}[dd] \ar@/_1pc/[lllddd]^{\epsilon_X} & \\
K_{\varphi}D_{\varphi}X \ar^-{\Psi_X}[rr] \ar^-{Leq_X}[rrru]
\ar^{\widehat{\epsilon}_X}[dd]  &  & EX \ar_-{eq'_X}[ru]
\ar@{-->}^{\nu_X}[dd]
 & & G^2X \\
  &  & &  GX \ar@<-.5ex>_{\Delta_X}[ru] \ar@<.5ex>^{Gx}[ru] \ar@{}^(.7){(1)}[uur]& \\
  X \ar@{=}[rr]   & & X \ar^{x}[ur] \ar@{}[lluu]|(.4){(2)}& & }
\end{equation}
En este diagrama, $\Psi$ est\'{a} determinada por $Leq_X$ en
virtud de la propiedad universal del igualador \eqref{ecuprima}.
Para definir $\nu_X$, necesitamos comprobar que el cuadrado $(1)$
conmuta serialmente. Pero eso es una consecuencia sencilla de la
naturalidad de $\varphi$ y del hecho de ser un homomorfismo de
com\'{o}nadas. De esta manera, $\nu_X$ est\'{a} determinado por
$\varphi_X \circ eq'_Y$ por la propiedad universal del igualador
can\'{o}nicamente asociado a la co\'{a}lgebra $(X,x)$ (ver
\cite[Proposition 3.3.4]{Barr/Wells:2002}). Comprobemos ahora que
el cuadrado $(2)$ es conmutativo, esto es, que
\begin{equation}\label{natiguales}
\widehat{\epsilon}_X = \nu_X \circ \Psi_X
\end{equation}
Es suficiente con comprobar que la igualdad \eqref{natiguales} es
cierta tras componer con el monomorfismo $x$. Realizamos
seguidamente dicho c\'{a}lculo:
\begin{equation*}
\begin{array}{lclr}
x \circ \widehat{\epsilon}_X & = & x \circ \epsilon_X \circ
Leq_X & \\
& = & \epsilon_{GX} \circ LRx \circ Leq_X & (\epsilon \hbox{ es
natural}) \\
& = & \epsilon_{GX} \circ LR\alpha_X  \circ Leq_X & (eq_X \hbox{ iguala } (Rx,\alpha_X))\\
& = & \epsilon_{GX} \circ LR\varphi_X \circ L\eta_{RX} \circ
Leq_X & \\
& = & \varphi_X \circ \epsilon_{LRX} \circ L\eta_{RX} \circ
Leq_X & (\epsilon \hbox{ es natural}) \\
& = & \varphi_X \circ Leq_X & \\
& = & \varphi_X \circ eq'_X \circ \Psi_X & \\
& = & x \circ \nu_X \circ \Psi_X
\end{array}
\end{equation*}

\begin{lemma}\label{clave2}
Supongamos que $\varphi_X$ es un monomorfismo. En el diagrama
\eqref{clave}, $\widehat{\epsilon}_X$ es un isomorfismo si, y
s\'{o}lo si, $\Psi_X$ y $\nu_X$ son isomorfismos.
\end{lemma}
\begin{proof}
Vamos a demostrar que si $\widehat{\epsilon}_X$ es un isomorfismo,
entonces
\begin{equation*}
\xymatrix{K_{\varphi}D_{\varphi}X \ar^{Leq_X}[rr] & & LRX
\ar@<.5ex>^-{L\alpha_X}[rr] \ar@<-.5ex>_-{LRx}[rr] & & LGRX}
\end{equation*}
es un igualador. Esto implica que $\Psi_X$ es entonces un
isomorfismo, de donde el lema se sigue f\'{a}cilmente. Sea, por tanto,
$f : Y \rightarrow LRY$ un homomorfismo en $\cat{A}$ tal que $LRx
\circ f = L\alpha_X \circ f$. Existe un \'{u}nico $\widehat{f} : Y
\rightarrow EX$ tal que $eq'_X \circ \widehat{f} = f$. Definimos
\begin{equation*}
\widetilde{f} = {\widehat{\epsilon}_X}^{\;-1} \circ \nu_X \circ
\widehat{f} : Y \longrightarrow K_{\varphi}D_{\varphi}X
\end{equation*}
Entonces
\begin{equation*}
\begin{array}{lclr}
\varphi_X \circ Leq_X \circ \widetilde{f} & = & \varphi_X \circ
Leq_X \circ
{\widehat{\epsilon}_X}^{\;-1} \circ \nu_X \circ \widehat{f} & \\
& = & \varphi_X \circ eq'_X \circ \Psi_X \circ
{\widehat{\epsilon}_X}^{\;-1} \circ
\nu_X \circ \widehat{f} & \\
& = & x \circ \nu_X \circ \Psi_X \circ \widehat{\epsilon}_X^{-1} \circ \nu_X \circ \widehat{f} \\
& = & x \circ \nu_X \circ \widehat{f} & (\hbox{por \eqref{natiguales}}) \\
& = & \varphi_X \circ eq'_X \circ \widehat{f}\\
& = & \varphi_X \circ f &
\end{array}
\end{equation*}
Como $\varphi_X$ es un monomorfismo, deducimos que $Leq_X \circ
\widetilde{f} = f$.
 Supongamos ahora que  $g : Y \rightarrow
K_{\varphi}D_{\varphi}X$ es tal que $Leq_X \circ g = f$. Basta con
demostrar que $\widehat{\epsilon}_X \circ g = \nu_X \circ
\widehat{f}$. Dado que $x$ es un monomorfismo, es suficiente con
que comprobemos que $x \circ \widehat{\epsilon}_X \circ g = x
\circ \nu_X \circ \widehat{f}$:
\begin{equation*}
\begin{array}{lclr} x \circ \widehat{\epsilon}_X \circ g & = & x \circ \nu_X
\circ
\Psi_X \circ g & (\hbox{por \eqref{natiguales}}) \\
& = & \varphi_X \circ eq'_X \circ \Psi_X \circ g & \\
& = & \varphi_X \circ Leq_X \circ g & \\
& = & \varphi_X \circ f \\
& = & \varphi_X \circ eq'_X \circ \widehat{f} & \\
& = & x \circ \nu_X \circ \widehat{f} &
\end{array}
\end{equation*}
\end{proof}

\begin{remark}
Obviamente, el Lema \ref{clave2} puede re-enunciarse diciendo que,
bajo la hip{\'o}tesis de ser $\varphi_X$ un monomorfismo,
$\widehat{\epsilon}_X$ es un isomorfismo si, y s\'{o}lo si, $L$
preserva el igualador \eqref{Dphi}, y $\nu_X$ es un isomorfismo.
\end{remark}

Recordemos que un funtor adjunto por la derecha es fiel y pleno si
y s\'{o}lo si la counidad es un isomorfismo \cite[Proposition
3.4.1]{Borceux:1994}

\begin{theorem}\label{Descent}
Sea $L : \cat{B} \rightarrow \cat{A}$ un funtor con un adjunto por
la derecha $R : \cat{A} \rightarrow \cat{B}$, y $G : \cat{A}
\rightarrow \cat{A}$ cualquier com\'{o}nada. Consideremos un funtor
$K_{\varphi} : \cat{B} \rightarrow \coalg{\cat{A}}{G}$ que hace
conmutar \eqref{UKL} con homomorfismo de com\'{o}nadas correspondiente
$\varphi : LR \rightarrow G$, y sean $\alpha : R \rightarrow RG$,
$\beta : L \rightarrow GL$ sus representaciones. Supongamos que
cada $G$--co\'{a}lgebra $(X,x)$, existe en $\cat{B}$ el igualador de
$\alpha_X, Rx$. Entonces el funtor adjunto por la izquierda
$D_{\varphi} : \coalg{\cat{A}}{G} \rightarrow \cat{B}$ al funtor
$K_{\varphi}$ es fiel y pleno si, y s\'{o}lo si, $L$ preserva los
igualadores de la forma \eqref{Dphi} y $\varphi$ es un isomorfismo
de com\'{o}nadas.
\end{theorem}
\begin{proof}
Sea $X$ cualquier objeto de $\cat{A}$. Un sencillo c\'{a}lculo
muestra que
\begin{equation}\label{ecucontr}
\xymatrix{RX \ar^-{\alpha_X}[rr] & &
\ar@/^1pc/^-{R\varepsilon_X}[ll] RGX \ar@<.5ex>^-{\alpha_{GX}}[rr]
\ar@<-.5ex>_-{R\Delta_X}[rr] & & RG^2X
\ar@/^1.5pc/^-{RG\varepsilon_X}[ll]}
\end{equation}
es un igualador contractible en el sentido de \cite[Section
3.3]{Barr/Wells:2002}. En efecto,
\begin{equation*}
R\varepsilon_X \circ \alpha_X = id_{RX} \quad (\hbox{por
\eqref{comod}}),
\end{equation*}
\begin{equation*}
RG\varepsilon_X \circ R \Delta_X  = id_{RGX} \quad (\hbox{por
\eqref{comod}}),
\end{equation*}
\begin{equation*}
RG\varepsilon_X \circ \alpha_{GX}  =  \alpha_X \circ
R\varepsilon_X \quad (\alpha \hbox{ es natural}).
\end{equation*}
El igualador \eqref{ecucontr} muestra que $eq_{GX} = \alpha_X$ y
$D_{\varphi}GX = RX$. De aqu\'{\i}, aplicando el funtor $L$ a
\eqref{ecucontr}, y en vista de \eqref{counidad1}, obtenemos que
$\widehat{\epsilon}_{GX} = \epsilon_{GX} \circ L\alpha_X =
\varphi_X$. Por tanto, si $D_{\varphi}$ es fiel y pleno, entonces
$\widehat{\epsilon}_{GX} = \varphi_X$ es un isomorfismo para toda
$G$--co\'{a}lgebra $(X,x)$. Del Lema \ref{clave2} deducimos
tambi\'{e}n que $L$ preserva el igualador \eqref{Dphi}. El
rec\'{\i}proco es consecuencia directa del Lema \ref{clave2} y la
ecuaci\'{o}n \eqref{natiguales}.
\end{proof}

Un isomorfismo de com\'{o}nadas $\varphi : LR \rightarrow G$ induce
\cite[Theorem 3.3]{Barr/Wells:2002} una equivalencia de categor\'{\i}as
$\coalg{\cat{A}}{LR} \cong \coalg{\cat{A}}{G}$. Combinando este
hecho con el Teorema \ref{Descent} y el Teorema de Beck
\cite[Theorem 3.10]{Barr/Wells:2002}, podemos obtener el Teorema
\ref{comonadic}. Preferimos, si el lector nos lo permite, incluir
una demostraci\'{o}n expl\'{\i}cita.

\begin{theorem}\label{comonadic}
Sea $L : \cat{B} \rightarrow \cat{A}$ un funtor con un adjunto por
la derecha $R : \cat{A} \rightarrow \cat{B}$, y $G : \cat{A}
\rightarrow \cat{A}$ cualquier com\'{o}nada. Consideremos un funtor
$K_{\varphi} : \cat{B} \rightarrow \coalg{\cat{A}}{G}$ que hace
conmutar \eqref{UKL} con homomorfismo de com\'{o}nadas correspondiente
$\varphi : LR \rightarrow G$, y sean $\alpha : R \rightarrow RG$,
$\beta : L \rightarrow GL$ sus representaciones. Supongamos que
cada $G$--co\'{a}lgebra $(X,x)$, existe en $\cat{B}$ el igualador de
$\alpha_X, Rx$. Entonces el funtor $K_{\varphi}$ es una
equivalencia de categor\'{\i}as entre $\cat{B}$ y $\coalg{\cat{A}}{G}$
si, y s\'{o}lo si, $L$ preserva los igualadores de la forma
\eqref{Dphi}, refleja isomorfismos, y $\varphi$ es un isomorfismo
de com\'{o}nadas.
\end{theorem}
\begin{proof}
Observemos primero que para cada objeto $Y$ de $\cat{B}$, la
unidad $\widehat{\eta}_Y$ de la adjunci\'{o}n $K_{\varphi} \dashv
D_{\varphi}$ demostrada en la Proposici\'{o}n \ref{adjuncion} viene
dada, seg\'{u}n la Observaci\'{o}n \ref{unitcounit}, como el igualador en
la fila horizontal del diagrama
\begin{equation}\label{unitphi2}
\xymatrix{ & & & RLRLY \ar^-{R\varphi_{LY}}[dr] & \\
D_{\varphi}K_{\varphi}Y \ar[rr] & & RLY
\ar@<.5ex>^-{\alpha_{LY}}[rr] \ar@<-.5ex>_-{R\beta_Y}[rr]
\ar@<.5ex>^-{\eta_{RLY}}[ur] \ar@<-.5ex>_-{RL\eta_Y}[ur]& &
RGLY\\
& Y \ar^-{\widehat{\eta}_Y}[ul] \ar_{\eta_Y}[ur] & & &}
\end{equation}
Si aplicamos el funtor $L$ al diagrama conmutativo
\eqref{unitphi2} obtenemos el diagrama
\begin{equation}\label{Lunitphi2}
\xymatrix{ & & & LRLRLY \ar^-{LR\varphi_{LY}}[dr] \ar@/_2.5pc/_{\epsilon_{LRLY}}[ld]& \\
LD_{\varphi}K_{\varphi}Y \ar[rr] & & LRLY
\ar@/_1pc/_-{\epsilon_{LY}}[ld] \ar@<.5ex>^-{L\alpha_{LY}}[rr]
\ar@<-.5ex>_-{LR\beta_Y}[rr] \ar@<.5ex>^-{L\eta_{RLY}}[ur]
\ar@<-.5ex>_-{LRL\eta_Y}[ur] & &
LRGLY\\
& LY \ar^-{L\widehat{\eta}_Y}[ul] \ar_{L\eta_Y}[ur], & & &}
\end{equation}
tambi\'{e}n conmutativo. Aqu\'{\i}, los morfismos $\epsilon_{LRLY}$ y
$\epsilon_{LY}$ hacen que la diagonal sea un igualador
contractible. Si $K_{\varphi}$ es una equivalencia de categor\'{\i}as,
entonces su adjunto por la derecha $D_{\varphi}$ es obviamente
fiel y pleno y, en virtud del Teorema \ref{Descent}, $\varphi$ es
un isomorfismo natural y $L$ preserva los igualadores de la forma
\eqref{Dphi}. Puesto que el funtor que olvida $U_G :
\coalg{\cat{A}}{G} \rightarrow \cat{A}$ refleja isomorfismos
\cite[Proposition 3.3.1]{Barr/Wells:2002}, deducimos de $L = U_G
\circ K_{\varphi}$ que $L$ refleja isomorfismos. Rec\'{\i}procamente,
si $\varphi$ es un isomorfismo natural y $L$ preserva los
igualadores \eqref{Dphi} y refleja isomorfismos, entonces, por el
Teorema \ref{Descent}, la counidad de la adjunci\'{o}n $K_{\varphi}
\dashv D_{\varphi}$ es un isomorfismo. Del diagrama
\eqref{Lunitphi2} deducimos que $L\widehat{\eta}_Y$ es un
isomorfismo y, como $L$ refleja isomorfismos, $\widehat{\eta}_Y$
ha de ser un isomorfismo, con lo que hemos demostrado que la
unidad de la adjunci\'{o}n $K_{\varphi} \dashv D_{\varphi}$ es tambi\'{e}n
un isomorfismo natural. Por tanto, $K_{\varphi}$ es una
equivalencia de categor\'{\i}as.
\end{proof}

\section{Coanillos sobre anillos firmes}

Sea $A$ un anillo, del que no suponemos posea un uno. Por
$\rmodu{A}$ denotamos la categor\'{\i}a de todos los $A$--m\'{o}dulos por
la derecha. Podemos tambi\'{e}n considerar m\'{o}dulos por la izquierda o
bim\'{o}dulos. El producto tensor sobre $A$ se denotar\'{a} mediante $-
\tensor{A} - $. Un $A$--m\'{o}dulo por la derecha $M$ se dir\'{a}
\emph{firme} \cite{Quillen:unp1997} si el homomorfismo
<<multiplicaci\'{o}n>> $\varpi_M^+ : M \tensor{A} A \rightarrow M$ es
biyectivo. Su inverso se denotar\'{a} por $d_M^+ : M \rightarrow M
\tensor{A} A$. Los m\'{o}dulos firmes por la izquierda se definen
an\'{a}logamente, con notaciones $\varpi_M^-$ y $d_M^-$ para el
homomorfismo <<multiplicaci\'{o}n>> y su inverso, respectivamente.
Obviamente, $\varpi_A^+ = \varpi_A^-$, con lo que, en caso de ser
este homomorfismo multiplicaci\'{o}n biyectivo, se tiene $d_A^+ =
d_A^-$. Diremos entonces que $A$ es un \emph{anillo firme}. Si $A$
es firme, entonces la subcategor\'{\i}a plena $\rmod{A}$ de $\rmodu{A}$
cuyos objetos son todos los m\'{o}dulos firmes una categor\'{\i}a abeliana
\cite[(4.6)]{Quillen:unp1997}, \cite[Corollary
1.3]{Grandjean/Vitale:1998}, \cite[Proposition 2.7]{Marin:1998},
aunque en general los igualadores no necesariamente se calculan en
grupos abelianos, debido esencialmente a la falta de exactitud del
funtor $ - \tensor{A} A$.

Supongamos que $A$ es un anillo firme. Un \emph{$A$--coanillo}
\cite{Sweedler:1975} es un $A$--bim\'{o}dulo $\coring{C}$ dotado de
dos homomorfismos de $A$--bim\'{o}dulos $\Delta_{\coring{C}} :
\coring{C} \rightarrow \coring{C} \tensor{A} \coring{C}$, y
$\varepsilon_{\coring{C}} : \coring{C} \rightarrow A$ que
verifican las siguientes ecuaciones:
\begin{equation}\label{coasociativa}
(\coring{C} \tensor{A} \Delta_{\coring{C}}) \circ
\Delta_{\coring{C}} = (\Delta_{\coring{C}} \tensor{A} \coring{C})
\circ \Delta_{\coring{C}}
\end{equation}
\begin{equation}\label{counitaria}
(\varepsilon_{\coring{C}} \tensor{A} \coring{C}) \circ
\Delta_{\coring{C}} = d_A^-, \qquad (\coring{C} \tensor{A}
\varepsilon_{\coring{C}}) \circ \Delta_{\coring{C}} = d_A^+
\end{equation}
En \eqref{coasociativa} hemos considerado el isomorfismo natural
$\coring{C} \tensor{A} (\coring{C} \tensor{A} \coring{C}) \cong
(\coring{C} \tensor{A} \coring{C}) \tensor{A} \coring{C}$ como una
igualdad, denotando el <<valor com\'{u}n>> como $\coring{C} \tensor{A}
\coring{C} \tensor{A} \coring{C}$. Un c\'{a}lculo directo, usando
\eqref{coasociativa} y \eqref{counitaria} demuestra que cada
$A$--coanillo
$(\coring{C},\Delta_{\coring{C}},\varepsilon_{\coring{C}})$
determina una com\'{o}nada sobre $\rmod{A}$ definida por el funtor
\begin{equation*}
\xymatrix{ - \tensor{A} \coring{C} : \rmod{A} \ar[r] & \rmod{A}}
\end{equation*}
y las transformaciones naturales
\begin{equation*}
\xymatrix{- \tensor{A} \coring{C} \ar^-{-\tensor{A}
\Delta_{\coring{C}}} [rr] & & - \tensor{A} \coring{C} \tensor{A}
\coring{C} \\
- \tensor{A} \coring{C} \ar^-{- \tensor{A}
\varepsilon_{\coring{C}}}[rr] & & - \tensor{A} A \cong id }
\end{equation*}
La categor\'{\i}a de co\'{a}lgebras para esta com\'{o}nada no es sino la
categor\'{\i}a $\rcomod{\coring{C}}$ de los $\coring{C}$--com\'{o}dulos por
la derecha.

Supongamos ahora dado un segundo anillo firme $B$, y un
$B-A$--bim\'{o}dulo $\Sigma$ firme.  Argumentando como en
\cite{Gomez/Vercruysse:2005arXiv}, tenemos un par de funtores
adjuntos
\begin{equation}\label{adjcoring}
\xymatrix{\rmod{B} \ar@<0.5ex>^-{\dostensor{-}{B}{\Sigma}}[rrr] &
& & \rmod{A},
\ar@<0.5ex>^-{\dostensor{\hom{A}{\Sigma}{-}}{B}{B}}[lll] &
\dostensor{-}{B}{\Sigma} \dashv
\dostensor{\hom{A}{\Sigma}{-}}{B}{B}}
\end{equation}
Ser\'{a} \'{u}til hacer expl\'{\i}citas la unidad y counidad de la adjunci\'{o}n
\eqref{adjcoring}. Para ello, dado $y \in Y$ para $Y$ un
$B$--m\'{o}dulo por la derecha firme, utilizaremos la notaci\'{o}n
$d_Y^+(y) = y^b \tensor{B} b \in Y \tensor{B} B$ (suma
sobreentendida). Por supuesto, este elemento del producto tensor
est\'{a} determinado por la condici\'{o}n $y^bb = y$. La counidad de la
adjunci\'{o}n es
\begin{equation}\label{unit3}
\xymatrix{\eta_Y : Y \ar[r] &
\dostensor{\hom{A}{\Sigma}{\dostensor{Y}{B}{\Sigma}}}{B}{B}, &
\eta_Y(y) = \dostensor{(\dostensor{y^b}{B}{-})}{B}{b},}
\end{equation}
y la counidad
\begin{equation}\label{counit3}
\xymatrix{\epsilon_X :
\trestensor{\hom{A}{\Sigma}{X}}{B}{B}{B}{\Sigma} \ar[r] & X, &
\epsilon_X(\trestensor{f}{B}{b}{B}{u}) = f(bu)}.
\end{equation}
Tenemos entonces la com\'{o}nada asociada
\begin{equation*}
(\trestensor{\hom{A}{\Sigma}{-}}{B}{B}{B}{\Sigma},
\dostensor{\eta_{\dostensor{\hom{A}{\Sigma}{-}}{B}{B}}}{B}{\Sigma},\epsilon)
\end{equation*}
Supongamos ahora dada una estructura de $B -
\coring{C}$--bicom\'{o}dulo sobre $\Sigma$,  esto es, un homomorfismo
de $B-A$--bim\'{o}dulos $\varrho_\Sigma : \Sigma \rightarrow \Sigma
\tensor{A} \coring{C}$ que verifica
\begin{equation}\label{rcomod}
(\varrho_{\Sigma} \tensor{A} \coring{C}) \circ \varrho_{\Sigma} =
(\Sigma \tensor{A} \Delta_{\coring{C}}) \circ \varrho_{\Sigma},
\qquad (\Sigma \tensor{A} \varepsilon_{\coring{C}})\circ
\varrho_{\Sigma} = d_{\Sigma}^+
\end{equation}
Se sigue f\'{a}cilmente de \eqref{rcomod} que la transformaci\'{o}n
natural
\begin{equation}\label{(C)}
\xymatrix{\beta : - \tensor{B} \Sigma \ar^-{- \tensor{B}
\varrho_{\Sigma}}[rr] & & - \tensor{B} \Sigma \tensor{A}
\coring{C}}
\end{equation}
satisface las condiciones de la afirmaci\'{o}n \emph{(C)} de la
Proposici\'{o}n \ref{correspondencia} y, en virtud de dicho Teorema,
da lugar a un homomorfismo can\'{o}nico de com\'{o}nadas $\mathsf{can}$
definido como la composici\'{o}n
\begin{equation*}
\xymatrix{\trestensor{\hom{A}{\Sigma}{-}}{B}{B}{B}{\Sigma}
\ar^-{\trestensor{\hom{A}{\Sigma}{-}}{B}{B}{B}{\varrho_{\Sigma}}}
[rrr] \ar_-{\mathsf{can}}[rrrd] & & &
\fourtensor{\hom{A}{\Sigma}{-}}{B}{B}{B}{\Sigma}{A}{\coring{C}}
\ar^{\epsilon \tensor{A} \coring{C}}[d]
\\ & & & - \tensor{A} \coring{C},}
\end{equation*}
o, si usamos una notaci\'{o}n de Heynemann-Sweedler abreviada,
tenemos, para cada $B$--m\'{o}dulo por la derecha $X$:
\begin{equation*}
\xymatrix{\trestensor{\hom{A}{\Sigma}{X}}{B}{B}{B}{\Sigma}
\ar^-{\mathsf{can}_X}[rr] & & X \tensor{A} \coring{C}
\\
\trestensor{f}{B}{b}{B}{u} \ar@{|->}[rr] & &
f(bu_{[0]})\tensor{A}u_{[1]} }
\end{equation*}
donde $\varrho_{\Sigma}(u) = u_{[0]} \tensor{A} u_{[1]}$ (suma
sobreentendida).

\medskip

 Podemos ahora aplicar la Proposici\'{o}n \ref{adjuncion} y el Teorema
\ref{Descent} para obtener

\begin{theorem}
El funtor $ - \tensor{B} \Sigma : \rmod{B} \rightarrow
\rcomod{\coring{C}}$ tiene un adjunto por la derecha
\begin{equation*}
\xymatrix{\dostensor{\hom{\coring{C}}{\Sigma}{-}}{B}{B} :
\rcomod{\coring{C}} \ar[rr] & &  \rmod{B}}
\end{equation*}
Este funtor es fiel y pleno si, y s\'{o}lo si, $\mathsf{can}$ es
un isomorfismo natural y $- \tensor{B} \Sigma :\rmod{B}
\rightarrow \rmod{A}$ preserva el igualador
\begin{equation}\label{iguamodu}
\xymatrix{ \hom{\coring{C}}{\Sigma}{X} \tensor{B} B \ar[r]  &
\hom{A}{\Sigma}{X}\tensor{B} B \ar@<.5ex>^-{\alpha_X}[rrr]
\ar@<-.5ex>_-{\hom{A}{\Sigma}{\varrho_X}\tensor{B}{B}}[rrr] & & &
\hom{A}{\Sigma}{X \tensor{A} \coring{C}} \tensor{B} B}
\end{equation}
para cada $\coring{C}$--com\'{o}dulo por la derecha $(X,\varrho_X)$,
donde $\alpha_X(f \tensor{B} b) = [(f \tensor{A} \coring{C}) \circ
\varrho_{\Sigma}] \tensor{B} b$ para $f \tensor{B} b \in
\hom{A}{\Sigma}{X} \tensor{B} B$.
\end{theorem}
\begin{proof}
Basta con comprobar que, en la presente situaci\'{o}n, $\alpha_X$, tal
como aparece definido en general en \eqref{tabladef}, est\'{a} dado en
la manera que afirma el enunciado, y que $- \tensor{B} B$ es
exacto por la izquierda puesto que es adjunto por la derecha del
funtor inclusi\'{o}n $J : \rmod{B} \rightarrow \rmodu{B}$.
\end{proof}

El Teorema \ref{comonadic} tiene la siguiente consecuencia:

\begin{theorem}
Sea $\Sigma$ un $B-\coring{C}$--bicom\'{o}dulo. El funtor $ -
\tensor{B} \Sigma : \rmod{B} \rightarrow \rcomod{\coring{C}}$ es
una equivalencia de categor\'{\i}as si, y s\'{o}lo si, $-
\tensor{B} \Sigma : \rmod{B} \rightarrow \rmod{A}$ preserva los
igualadores de la forma \eqref{iguamodu}, refleja isomorfismos y
la transformaci{\'o}n can\'{o}nica $\mathsf{can}$ es un isomorfismo
natural.
\end{theorem}

Vamos ahora a deducir algunos de los resultados fundamentales de
\cite{Gomez/Vercruysse:2005arXiv}. Escribamos $\Sigma^* =
\hom{A}{\Sigma}{A}$, y supongamos que $A$ es un anillo con unidad.
El $B$--bim\'{o}dulo $\Sigma \tensor{A} \Sigma^*$ tiene una estructura
de $B$--anillo (sin uno, en general), con multiplicaci\'{o}n
asociativa definida por
\begin{equation*}
\mu(x \tensor{A} \phi \tensor{B} y \tensor{A} \psi) = x \phi(y)
\tensor{A} \psi = x \tensor{A} \phi(y) \psi, \qquad (x \tensor{A}
\phi \tensor{B} y \tensor{A} \psi \in \Sigma \tensor{A} \Sigma^*
\tensor{B} \Sigma \tensor{A} \Sigma^*)
\end{equation*}
 La situaci\'{o}n tratada en
\cite{Gomez/Vercruysse:2005arXiv} parte de un homomorfismo de
anillos $\iota : B \rightarrow \Sigma \tensor{A} \Sigma^*$, lo que
permite demostrar que existe un isomorfismo natural
\begin{equation}\label{isomor}
\hom{A}{\Sigma}{-} \tensor{B} B \simeq - \tensor{A} \Sigma^*
\tensor{B} B \qquad (h \tensor{B} b \mapsto h(e_c) \tensor{A}
e^*_c \tensor{B} b^c)
\end{equation}
Aqu\'{\i}, estamos usando la notaci\'{o}n $\iota (b) = e_b \tensor{B}
e^*_b$ para $b \in B$ (suma sobreentendida).
  De \eqref{isomor} deducimos, usando la notaci\'{o}n
$\Sigma^{\dag} = \Sigma^* \tensor{B} B$,  una adjunci\'{o}n
\begin{equation}\label{adj6}
\xymatrix{\rmod{B} \ar@<0.5ex>^-{\dostensor{-}{B}{\Sigma}}[rrr] &
& & \rmod{A}, \ar@<0.5ex>^-{\dostensor{-}{A}{\Sigma^{\dag}}}[lll]
& \dostensor{-}{B}{\Sigma} \dashv \dostensor{-}{A}{\Sigma^{\dag}}}
\end{equation}
cuya counidad es
\begin{equation}\label{counit6}
\xymatrix{\epsilon_M : \trestensor{M}{A}{\Sigma^{\dag}}{B}{\Sigma}
\ar[r] & M & (\fourtensor{m}{A}{\phi}{B}{b}{B}{x} \mapsto
m\phi(bx)),}
\end{equation}
y cuya unidad es
\[
\xymatrix{\eta_N : N \ar[r] &
\trestensor{N}{B}{\Sigma}{A}{\Sigma^{\dag}} & (n \mapsto
\fourtensor{n^b}{B}{e_c}{A}{e^*_c}{B}{b^c})}
\]
A partir del par adjunto \eqref{adj6} tenemos la com\'{o}nada sobre
$\rmod{A}$
\begin{equation}\label{comatrixcomon}
(- \tensor{A} \Sigma^{\dag} \tensor{B} \Sigma, \eta_{- \tensor{A}
\Sigma^{\dag}} \tensor{B} \Sigma, \epsilon )
\end{equation}
Evaluando en $A$, obtenemos el $A$--coanillo $A \tensor{A}
\Sigma^{\dag} \tensor{B} \Sigma \cong \Sigma^{\dag} \tensor{B}
\Sigma$ con comultiplicaci\'{o}n $\Delta^{\dag} = \eta_{A \tensor{A}
\Sigma^{\dag}} \tensor{B} \Sigma$ y counidad $\epsilon_A$.
Expl\'{\i}citamente,
\begin{equation*}
\Delta^{\dag}(\trestensor{\phi}{B}{b}{B}{u}) =
\fourtensor{\phi}{B}{c}{B}{\abrir{d}}{R}\dostensor{(b^c)^d}{B}{u}
\end{equation*}
\begin{equation*}
\epsilon_A(\phi \tensor{B} b \tensor{B} u) = \phi(bu)
\end{equation*}
Adem\'{a}s, la com\'{o}nada \eqref{comatrixcomon} est\'{a} determinada por el
\emph{coanillo de comatrices} $(\Sigma^{\dag} \tensor{B} \Sigma,
\Delta^{\dag},\epsilon_A)$, ya que el funtor subyacente a la
com\'{o}nada es un producto tensor sobre $A$.

Si volvemos ahora al caso en que $(\Sigma, \varrho_{\Sigma})$ es
un $B-\coring{C}$--bicom\'{o}dulo, tenemos que la
transformaci\'{o}n natural \eqref{(C)} da lugar, en virtud de la
Proposici\'{o}n \ref{correspondencia} a un homomorfismo de
com\'{o}nadas $\mathsf{can}^{\dag}$ definido (ver
\eqref{tabladef}) en $(X,\varrho_X)$ por
\begin{equation*}
{\mathsf{can}}^{\dag}_X = (X \tensor{A} \epsilon_A \tensor{A}
\coring{C}) \circ (X \tensor{A} \Sigma^{\dag} \tensor{B}
\varrho_{\Sigma})
\end{equation*}
Obviamente,
\begin{equation}\label{canplus}
\mathsf{can}^{\dag}_X = X \tensor{A} \mathsf{can}^{\dag}_A
\end{equation}
donde $\mathsf{can}^{\dag}_A$ es la aplicaci\'{o}n
\begin{equation}
\xymatrix{\mathsf{can}^{\dag}_A : \Sigma^{\dag} \tensor{B} \Sigma
\ar[r] & \coring{C}, & & \phi \tensor{B} b \tensor{B} u \mapsto
\phi(bu_{[0]})u_{[1]}}
\end{equation}
que es un homomorfismo de $A$--coanillos, puesto que
$\mathsf{can}^{\dag}$ es un homomorfismo de com\'{o}nadas. De
acuerdo con la Proposici\'{o}n \ref{adjuncion}, el funtor de
comparaci\'{o}n $- \tensor{B} \Sigma : \rmod{B} \rightarrow
\rcomod{\coring{C}}$ tiene un adjunto por la derecha, que, sobre
un $\coring{C}$--com\'{o}dulo por la derecha $(X, \varrho_X)$,
est\'{a} definida como el igualador en $\rmod{B}$ del par
$(\alpha_X,\varrho_X \tensor{A} \Sigma^{\dag})$. Un sencillo
c\'{a}lculo muestra que, en el presente caso, $\alpha_X = X
\tensor{A} \alpha_{\Sigma^{\dag}}$, donde
\begin{equation}\label{sigmaplus}
\alpha_{\Sigma^{\dag}} = (\mathsf{can}^{\dag}_A \tensor{A}
\Sigma^{\dag}) \circ \eta_{\Sigma^{\dag}}
\end{equation}
De los diagramas \eqref{comod} deducimos inmediatamente que
$(\Sigma^{\dag},\alpha_{\Sigma^{\dag}})$ es un
$\coring{C}-B$--bicom\'{o}dulo. De esta manera, el funtor definido en
la Proposici\'{o}n \ref{adjuncion} es un producto cotensor, ya que
est\'{a} definido por el igualador
\begin{equation}\label{cotensor}
\xymatrix{X \cotensor{\coring{C}} \Sigma^{\dag} \ar[r] & X
\tensor{A} \Sigma^{\dag} \ar@<.5ex>^-{X \tensor{A}
\alpha_{\Sigma^{\dag}}}[rr] \ar@<-.5ex>_-{\varrho_X \tensor{A}
\Sigma^{\dag}}[rr] & & X \tensor{A} \coring{C} \tensor{A}
\Sigma^{\dag} }
\end{equation}
De la Proposici\'{o}n \ref{adjuncion} y el Teorema \ref{Descent}
deducimos, pues:

\begin{theorem}\cite{Gomez/Vercruysse:2005arXiv} Sea $\coring{C}$
un $A$--coanillo y $\Sigma$ un $B-\coring{C}$--com\'{o}dulo, firme
como $B$--m\'{o}dulo por la izquierda. Si $\iota : B \rightarrow
\Sigma \tensor{A} \Sigma^*$ es un homomorfismo de $B$--anillos,
entonces el funtor $- \tensor{B} \Sigma : \rmod{B} \rightarrow
\rcomod{\coring{C}}$ tiene como adjunto por la derecha al funtor
producto cotensor $- \cotensor{\coring{C}} \Sigma^{\dag} :
\rcomod{\coring{C}} \rightarrow \rmod{B}$. Este funtor es fiel y
pleno si, y s\'{o}lo si, $\mathsf{can}^{\dag}_A : \Sigma^{\dag}
\tensor{B} \Sigma \rightarrow \coring{C}$ es un isomorfismo de
$A$--coanillos y $ - \tensor{B} \Sigma : \rmod{B} \rightarrow
\rmod{A}$ preserva los igualadores de la forma \eqref{cotensor}.
\end{theorem}

El Teorema \ref{comonadic} da en la presente situaci\'{o}n:

\begin{theorem} Sea $\coring{C}$
un $A$--coanillo y $\Sigma$ un $B-\coring{C}$--com\'{o}dulo, firme
como $B$--m\'{o}dulo por la izquierda. Si $\iota : B \rightarrow
\Sigma \tensor{A} \Sigma^*$ es un homomorfismo de $B$--anillos,
entonces el funtor $- \tensor{B} \Sigma : \rmod{B} \rightarrow
\rcomod{\coring{C}}$ es una equivalencia de categor\'{\i}as si, y
s\'{o}lo si, $\mathsf{can}^{\dag}_A : \Sigma^{\dag} \tensor{B}
\Sigma \rightarrow \coring{C}$ es un isomorfismo de $A$--coanillos
y $ - \tensor{B} \Sigma : \rmod{B} \rightarrow \rmod{A}$ refleja
isomorfismos y preserva los igualadores de la forma
\eqref{cotensor}.
\end{theorem}

\begin{remark}
Cuando tanto $A$ como $B$ son anillos con uno, y $\Sigma$ es
finitamente generado y proyectivo como $A$--m\'{o}dulo por la derecha,
entonces $\Sigma \tensor{A} \Sigma^* \cong \rend{A}{\Sigma}$ e
$\iota : B \rightarrow \rend{A}{\Sigma}$ no es sino el homomorfismo
que lleva cada $b \in B$ en el endomorfismo multiplicaci\'{o}n por $b$.
De hecho, podemos tomar $B = \rend{\coring{C}}{\Sigma}$. Esta es la
situaci\'{o}n considerada en \cite[Section 3]{ElKaoutit/Gomez:2003}.
Cuando $\Sigma = A$, cada estructura de $\coring{C}$--com\'{o}dulo est\'{a}
determinada por un elemento <<como de grupo>> $g \in \coring{C}$, y
tenemos la situaci\'{o}n estudiada en \cite[Section 5]{Brzezinski:2002}.
Cuando $\mathsf{can}_A$ es un isomorfismo, se dice en
\cite{ElKaoutit/Gomez:2003} o \cite{Brzezinski:2002},
respectivamente, que $\coring{C}$ es un coanillo de Galois ($\Sigma$
\'{o} $g$ se sobrentienden), o que $\Sigma$ es un $\coring{C}$--com\'{o}dulo
de Galois \cite{Brzezinski/Wisbauer:2003}. El coanillo $\Sigma^*
\tensor{B} \Sigma$ se suele llamar coanillo de comatrices.
\end{remark}

\begin{center}
\textbf{AGRADECIMIENTOS}
\end{center}
Quiero dar las gracias M. Mesablishvili por advertirme de c{\'o}mo
buena parte de los resultados de las secciones 1 y 2 son
esencialmente conocidos, proporcion{\'a}ndome adem{\'a}s las referencias
adecuadas. Tambi{\'e}n agradezco a C. Menini su atenta lectura de la
primera versi{\'o}n de esta nota, que la llev{\'o} a indicarme un error en
la demostraci{\'o}n del anterior enunciado del Lema \ref{clave2}, y
otras inexactitudes.

\providecommand{\bysame}{\leavevmode\hbox
to3em{\hrulefill}\thinspace}
\providecommand{\MR}{\relax\ifhmode\unskip\space\fi MR }
\providecommand{\MRhref}[2]{%
  \href{http://www.ams.org/mathscinet-getitem?mr=#1}{#2}
} \providecommand{\href}[2]{#2}

\end{document}